\documentclass[12TP,draft]{article}

\begin{document}

\begin{center}
\LARGE\noindent\textbf{On a problem of Wang concerning the Hamiltonicity of bipartite digraphs}\\

\end{center}
\begin{center}
\noindent\textbf{Samvel Kh. Darbinyan and Iskandar A. Karapetyan }\\

\end{center}
\begin{center}

Institute for Informatics and Automation Problems,

 Armenian National Academy of Sciences

E-mail: samdarbin@ipia.sci.am, isko@ipia.sci.am\\
\end{center}

\textbf{Abstract}

 R. Wang (Discrete Mathematics and Theoretical Computer Science, vol. 19(3), 2017) proposed the following problem.
 
\textbf{Problem.} Let $D$ be a strongly connected balanced bipartite directed graph of order $2a\geq 8$. Suppose that 
$d(x)\geq 2a-k$, $ d(y)\geq a+k$ or $d(y)\geq 2a-k$, $ d(x)\geq a+k$ for every  pair of vertices $x,y$ with a common out-neighbour, where $2 \leq k\leq a/2$. Is  $D$  Hamiltonian?

 In this paper, we prove that if a digraph $D$ satisfies the conditions of this problem, then 

(i) $D$ contains a cycle factor, 

 (ii) for every vertex 
$x\in V(D)$ there exists a vertex $y\in V(D)$ such that $x$ and $y$ have a common out-neighbour.\\

\textbf{Keywords:} Digraph, cycle, Hamiltonian cycle, bipartite  digraph, perfect matching.\\

\section {Introduction} 

In this paper, we consider finite directed graphs (digraphs) without loops and multiple arcs.  A digraph $D$ is called  Hamiltonian  if it contains a  Hamiltonian cycle, i.e., a cycle that includes every vertex
 of $D$. The vertex set and the arc set of a digraph $D$ are  denoted
  by $V(D)$  and   $A(D)$, respectively. The order of a digraph $D$ is the number of its vertices. A cycle factor in $D$ is a collection of vertex-disjoint cycles $C_1, C_2, \ldots , C_l$ such that $V(C_1)\cup V(C_2)\cup \ldots \cup V(C_l)= V(D)$. 
A digraph $D$ is bipartite  if there exists a partition $X$, $Y$ of $V(D)$ into two partite sets such that every arc of $D$ has its end-vertices in different partite sets. 
It is called balanced if $|X|=|Y|$.

  There are a number of conditions that guarantee that a bipartite digraph is  Hamiltonian  (see,   e.g., \cite{[1]}-\cite{[11]}). Let us recall the following degree conditions that guarantee that a  balanced bipartite digraph is  Hamiltonian. \\

\noindent\textbf{Theorem 1.1} ( Adamus, Adamus and Yeo \cite{[8]}). 
 {\it Let $D$ be a balanced bipartite digraph of order $2a$, where $a\geq 2$. Then $D$ is Hamiltonian provided one of the following holds}:

(a) {\it $d(u)+d(v)\geq 3a+1$ for every pair  of non-adjacent  distinct vertices $u$ and $v$ of $D$};

(b) {\it $D$ is strongly connected and  $d(u)+d(v)\geq 3a$ for every pair 
 of non-adjacent  distinct vertices $u$ and $v$ of $D$};

(c) {\it the minimal degree of $D$ is at least $(3a+1)/2$};

(d) {\it $D$ is strongly connected and the minimal degree of $D$ is at least $3a/2$}.\\

Observe that   Theorem 1.1 imposes a degree condition on all pairs of non-adjacent vertices. In the following  theorems  a degree condition requires only  for some pairs of non-adjacent vertices. \\

 \textbf{Theorem 1.2} (J. Adamus \cite{[9]}).  {\it Let $D$ be a strongly connected balanced bipartite digraph of order $2a\geq 6$. If  $d(x)+ d(y)\geq 3a$ for every pair of vertices $x$, $y$ with a common out-neighbour or a common in-neighbour, then $D$ is Hamiltonian.}\\

 Notice that Theorem 1.2 improves Theorem 1.1.\\

Some sufficient conditions for the existence of  Hamiltonian cycles in a bipartite tournament are described in the survey paper \cite{[3]} by Gutin.
A characterization for hamiltonicity for semicomplete bipartite digraphs was obtained independently  by Gutin \cite{[2]}  and  H\"{a}ggkvist and Manoussakis \cite{[4]}.\\

\textbf{Theorem 1.3} (Wang \cite{[10]}). {\it Let $D$ be a strongly connected balanced bipartite digraph of order $2a$, where $a\geq 1$. Suppose that, for every  pair of vertices $\{x,y\}$ with a common out-neighbour, either $d(x)\geq 2a-1$ and $d(y)\geq a+1$ or $d(y)\geq 2a-1$ and $d(x)\geq a+1$. Then $D$ is Hamiltonian.}\\

Before stating the next theorem we need to define a balanced bipartite digraph of order eight.\\

 \textbf{Example 1.} Let $D(8)$ be a  bipartite digraph   with partite sets $X=\{x_0,x_1,x_2,x_3\}$ and 
$Y=\{y_0,y_1,y_2,y_3\}$, and the arc set  $A(D(8))$ contains exactly the following  arcs: $y_0x_1$, $y_1x_0$, $x_2y_3$, $x_3y_2$ and  all the arcs of the following 2-cycles: 
$x_i\leftrightarrow y_i$, $i\in [0,3]$, $y_0\leftrightarrow x_2$, $y_0\leftrightarrow x_3$, $y_1\leftrightarrow x_2$ and  $y_1\leftrightarrow x_3$. 

It is easy to see that 
$$
d(x_2)=d(x_3)=d(y_0)=d(y_1)=7 \quad \hbox{and} \quad d(x_0)=d(x_1)=d(y_2)=d(y_3)=3,
$$
and  the dominating pairs in $D(8)$ are: $\{y_0,y_1\}$, $\{y_0,y_2\}$,$\{y_0,y_3\}$,$\{y_1,y_2\}$, $\{y_1,y_3\}$, $\{x_0,x_2\}$,
$\{x_0,x_3\}$,
$\{x_1,x_2\}$, $\{x_1,x_3\}$ and $\{x_2,x_3\}$. Note that  $max \{d(x), d(y)\}\geq 2a-1$ for every dominating pair of vertices $x$, $y$. Since $x_0y_0x_3y_2x_2$ $y_1x_0$ is a cycle in $D(8)$, 
it is not difficult to check that $D(8)$ is strong. 

Observe that $D(8)$ is not Hamiltonian. Indeed, if $C$ is a Hamiltonian cycle in $D(8)$, then $C$ would  contain the arcs $x_1y_1$ and $x_0y_0$. Therefore, $C$ would  contain the path $x_1y_1x_0y_0$ or the path $x_0y_0x_1y_1$, which is impossible since $N^-(x_0)=N^-(x_1)=\{y_0,y_1\}$.

Notice that the digraph $D(8)$ does not satisfy the conditions of Wang's theorem.\\

 \textbf{Theorem 1.4} (Darbinyan \cite{[11]}). {\it Let $D$ be a strongly connected balanced bipartite digraph of order $2a\geq 8$. Suppose that  $max \{d(x), d(y)\}\geq 2a-1$ for every pair of vertices $x$, $y$ with a common out-neighbour.
 Then  $D$ is Hamiltonian unless $D$ is isomorphic to the digraph $D(8)$} ({\it for definition of $D(8)$, see Example 1}).\\

For $a\geq 4$ Theorem 1.4 improves Wang's theorem.\\

A digraph $D$ of order $n$ is called pancyclic if it contains cycles of every length $k$, $3\leq k\leq n$.
A balanced bipartite digraph of order $2a$ is called even pancyclic if it contains  cycles of every length $2k$,  $2\leq k \leq a$.

 There are various sufficient conditions for a digraph (undirected graph) to be  Hamiltonian  are also sufficient for the digraph (undirected graph) to be pancyclic \cite{[1]}.

Recently, the following results were proved.\\
 
\textbf{Theorem 1.5} (Darbinyan \cite{[12]}). {\it Let $D$ be a strongly connected balanced bipartite digraph of order $2a\geq 8$ other than a directed cycle of length $2a$. If $ max\{d(x), d(y)\}\geq 2a-1$ for every dominating pair of vertices $x,y$,  then either $D$ contains  cycles of all even lengths less than or equal to $2a$ or $D$ is isomorphic to the digraph $D(8)$. }\\

\textbf{Theorem 1.6} (Meszka \cite{[13]}). {\it Let $D$ be a balanced bipartite digraph of order $2a\geq 4$ with partite sets $X$ and $Y$. If $d(x)+d(y)\geq 3a+1$  for every  pair of distinct vertices $x,y$ either both in $X$ or both in $Y$,  then $D$ contains  cycles of all even lengths less than or equal to $2a$. }\\

\textbf{Theorem 1.7} (Darbinyan \cite{[14]}). {\it Let $D$ be a strongly connected balanced bipartite digraph of order $2a\geq 6$ with partite sets $X$ and $Y$. If $d(x)+d(y)\geq 3a$  for every  pair of distinct vertices $x,y$ either both in $X$ or both in $Y$,  then $D$ contains  cycles of all even lengths less than or equal to $2a$. }\\

\textbf{Theorem 1.8} (Adamus \cite{[15]}). {\it Let $D$ be a strongly connected balanced bipartite digraph of order $2a\geq 6$.  If $d(x)+d(y)\geq 3a$  for every  pair of distinct vertices $x,y$ with a common in-neighbour or a common out-neighbour,  then $D$ contains  cycles of all even lengths less than or equal to $2a$ or else $D$ is a directed cycle of length $2a$. }\\

\textbf{Definition 1.} {\it Let $D$ be a balanced bipartite digraph of order $2a$, where $a\geq 2$. For any integer
 $k\geq 0$, we will say that
 $D$ satisfies condition $B_k$ when 

$$
d(x)\geq 2a-k, d(y)\geq a+k \quad \hbox{or} \quad d(x)\geq a+k, d(y)\geq 2a-k
$$ 
for any dominating pair of vertices $x,y$ in $D$.}\\

In \cite{[10]}, Wang proposed the following problem.\\

\textbf{Problem} (Wang \cite{[10]}). {\it Let $D$ be a strongly connected balanced bipartite digraph of order $2a\geq 8$ satisfying the condition $B_k$ with $2\leq k\leq a/2$. Is $D$ Hamiltonian?}\\

Before stating the next theorems we need to define a digraph of order ten.  \\

\textbf{Example 2.} Let $D(10)$ be a bipartite digraph with partite sets $X=\{x_0,x_1,x_2,x_3,x_4\}$ and 
$Y=\{y_0,y_1,y_2,y_3,y_4\}$ satisfying the following conditions: 
The induced subdigraph $\langle\{x_1,x_2,x_3,y_0,y_4\}\rangle$ is a complete bipartite digraph with partite sets  $\{x_1,x_2,x_3\}$ and $\{y_0,y_4\}$;  $\{x_1,x_2,x_3\}\rightarrow\{y_1,y_2,y_3\}$; $x_4\leftrightarrow y_4$; $x_0\leftrightarrow y_0$, $x_3\leftrightarrow y_1$ and $x_i\leftrightarrow y_{i+1}$ for all $i\in [1,3]$. $D(10)$ contains no other arcs.\\

It is easy to check that the digraph $D(10)$ is strongly connected  and  $ max\{d(x), d(y)\}\geq 2a-2$ for every dominating pair of vertices $x,y$, but the underlying undirected graph of $D(10)$ is not 2-connected and $D(10)$ has no cycle of length 8. (It follows from the facts that $d(x_0)=d(x_4)=2$ and $x_0$ ($x_4$) is on 2-cycle). Since $x_1y_1x_3y_3x_2y_2x_1$ is a cycle of length 6, $x_0\leftrightarrow y_0$ and $x_4\leftrightarrow y_4$, it is not difficult to check that any digraph obtained from $D(10)$ by adding a new arc the  one end-vertex of which is $x_0$ or $x_4$ contains a cycle of length eight. Moreover, if to $A(D)$ we add some new arcs of the type $y_ix_j$, where $i\in [1,3]$ and $j\in [1,3]$, then we always  obtain a digraph which contains a dominating pair of vertices, say $u,v$, such that  $ max\{d(u), d(v)\}\leq 2a-3$. \\

\textbf{Theorem 1.9} (\cite{[16]}, \cite{[17]}). {\it Let $D$ be a  balanced bipartite strongly connected digraph of order $2a\geq 10$ other than a directed cycle of length $2a$. Suppose that  $D$ 
 $ max\{d(x), d(y)\}\geq 2a-2$ for every dominating pair of vertices $x,y$ of $D$.  Then  $D$ contains  cycles of all lengths $2, 4, \ldots , 2a-2$ unless $D$ is  isomorphic to the digraph $D(10)$.}\\

Clearly, the existence of a cycle factor is a necessary condition for a digraph to be Hamiltonian. In this note we prove the following theorem.\\

\textbf{Theorem 1.10}. {\it Let $D$ be a strongly connected balanced bipartite digraph of order $2a\geq 8$ satisfying the condition $B_k$ with $2\leq k\leq a/2$. Then  $D$ contains a cycle factor.}

\section {Terminology and Notation}

Terminology  and notation not described below follow \cite{[1]}.
If $xy\in A(D)$, then  we  say that $x$ dominates $y$ or $y$ is an out-neighbour of $x$, and $x$ is an in-neighbour of $y$.

 Let $x,y$ be distinct vertices in a digraph $D$. The pair $\{x,y\}$ is called dominating if there is a vertex $z$ in $D$ such that  $xz\in A(D)$ and $yz\in A(D)$. In this case we say that $x$ is a partner of $y$ and  $y$ is a partner of $x$. 
 If $x\in V(D)$ and $A=\{x\}$ we sometimes will write $x$ instead of $\{x\}$. 
$A\rightarrow B$ means that every vertex of $A$ dominates every vertex of $B$. The notation $x\leftrightarrow y$ denotes that $xy\in A(D)$ and $yx\in A(D)$.

 Let $N^+(x)$, $N^-(x)$ denote the set of  out-neighbours, respectively the set  of in-neighbours of a vertex $x$ in a digraph $D$.  If $A\subseteq V(D)$, then $N^+(x,A)=A\cap N^+(x)$, $N^-(x,A)=A\cap N^-(x)$ and $N^+(A)=\cup_{x\in A}N^+(x)$, $N^-(A)=\cup_{x\in A}N^-(x)$. The out-degree of $x$ is $d^+(x)=|N^+(x)|$ and $d^-(x)=|N^-(x)|$ is the in-degree of $x$. Similarly, $d^+(x,A)=|N^+(x,A)|$ and $d^-(x,A)=|N^-(x,A)|$. The degree of the vertex $x$ in $D$ is defined as $d(x)=d^+(x)+d^-(x)$ (similarly, $d(x,A)=d^+(x,A)+d^-(x,A)$).

The path (respectively, the cycle) consisting of the distinct vertices $x_1,x_2,\ldots ,x_m$ ( $m\geq 2 $) and the arcs $x_ix_{i+1}$, $i\in [1,m-1]$  (respectively, $x_ix_{i+1}$, $i\in [1,m-1]$, and $x_mx_1$), is denoted by  $x_1x_2\cdots x_m$ (respectively, $x_1x_2\cdots x_mx_1$). 
We say that $x_1x_2\cdots x_m$ is a path from $x_1$ to $x_m$ or is an $(x_1,x_m)$-path. 
 Given a vertex $x$ of a directed path $P$ or a directed cycle $C$, we denote by $x^+$ (respectively, by $x^-$) the successor (respectively, the predecessor) of $x$ (on $P$ or $C$), and in case of ambiguity, we precise $P$ or $C$ as a subscript (that is $x^+_P$ \ldots).

A digraph $D$ is strongly connected (or, just, strong) if there exists an  $(x,y)$-path in $D$ for every ordered pair of    distinct vertices $x,y$ of $D$.
    Two distinct vertices $x$ and $y$ are adjacent if $xy\in A(D)$ or $yx\in A(D) $ (or both). 

Let $H$ be a non-trivial proper subset of vertices  of a digraph $D$. An $(x,y)$-path $P$ is an $H$-bypass if $|V(P)|\geq 3$, $x\not=y$ and $V(P)\cap H=\{x,y\}$.

Let $D$ be a balanced bipartite  digraph with partite sets  $X$ and $Y$. A matching from $X$ to $Y$ is an independent set of arcs with origin in $X$ and terminus in $Y$. (A  set of arcs with no common end-vertices is called independent). If $D$ is balanced, one says that such a matching is perfect if it consists of precisely $|X|$ arcs.\\
The underlying undirected graph of a digraph $D$ is denoted by $UG(D)$, it contains an edge $xy$ 
if $xy\in A(D)$ or $yx\in A(D)$ (or both).

\section {Main result}

Theorem 1.10 is the main result of this paper.\\

\textbf{Proof of theorem 1.10.} 

Let $D$ be a digraph satisfying the conditions of the theorem. Ore in \cite{[18]} (Section 8.6) has shown that a balanced bipartite digraph $D$ with partite sets $X$ and $Y$ has a cycle factor if and only if $D$ contains a perfect matching from $X$ to $Y$ and a perfect matching from $Y$ to $X$.

Therefore, by the well-known K\"{o}ning-Hall theorem (see, e.g., \cite{[19]}) to show that $D$ contains a perfect matching from $X$ to $Y$, it suffices to show that $|N^+(S)|\geq |S|$ for every set $S\subseteq X$. Let $S\subseteq X$. If $|S|=1$ or $|S|=a$, then $|N^+(S)|\geq |S|$ since $D$ is strongly connected. Assume that $2\leq |S|\leq a-1$. We claim that $|N^+(S)|\geq |S|$. Suppose that this is not the case, i.e., $|N^+(S)|\leq |S|-1\leq a-2$. From this and strongly connectedness of  $D$ it follows that there are two vertices $x,y\in S$ and a vertex $z\in N^+(S)$ such that $\{x,y\}\rightarrow z$, i.e., $\{x,y\}$ is a dominating pair. Therefore, by condition $B_k$,
$d(x)\geq 2a-k$ and $d(y)\geq a+k$ or $d(x)\geq a+k$ and $d(y)\geq 2a-k$. 
Without loss of generality, we assume that $d(x)\geq 2a-k$ and $d(y)\geq a+k$. Then 
$$
2a-k\leq d(x)\leq 2|N^+(S)|+a -|N^+(S)|=a+|N^+(S)|.
$$
Therefore, $|N^+(S)|\geq a-k$ and $|S|\geq a-k+1$. 

\textbf{Proposition 1.} Let $\{u,v\}$ be a dominating pair of vertices of $D$. Then from condition $B_k$ and 
$2\leq k\leq a/2$ it follows that $d(u)\geq a+k$ and $d(v)\geq a+k$, i.e., if a vertex $z$ has a partner in $D$ , then $d(z)\geq a+k$.\\

We claim that each vertex in $Y\setminus N^+(S)$ has no partner in $D$. Indeed, let $u$ be an arbitrary vertex in   $Y\setminus N^+(S)$. Since $|S|\geq a-k+1$, we have
$$
d(u)\leq |S|+2(a-|S|)=2a-|S|\leq a+k-1,
$$
which contradicts Proposition 1. This means that $u$ has no partner in $D$.\\

Without loss of generality,  assume that
$$
S= \{x_1,x_2, \ldots , x_s\} \quad  \hbox{and}  \quad N^+(S)=\{y_1,y_2, \ldots , y_t\}.
$$
Recall that every vertex $y_i$ with $t+1\leq i\leq a$ has no partner in $D$.
Note that $s\geq t+1$, $a-s\leq a-t-1$ and  there is no arc from a vertex of $\{x_1,x_2, \ldots , x_s\}$ to a vertex of 
$\{y_{t+1},y_{t+2}, \ldots , y_a\}$.
From this and strongly connectedness of $D$ it follows that there is a vertex $x_{i_1}$ such that $y_{t+1}x_{i_1}\in A(D)$.
Since $y_{t+1}$ has no partner, it follows that $d^-(x_{i_1}, Y\setminus \{y_{t+1}\})=0$. 
Therefore, $d(x_{i_1})\leq a+1\leq a+k-1$ since $k\geq 2$.
By Proposition 1, this means that the vertex $x_{i_1}$ also has no partner. Since $D$ is strongly connected, there is a vertex $y_{i_2}\in Y$
such that $x_{i_1}y_{i_2}\in A(D)$. Then $d^-(y_{i_2}, X\setminus \{x_{i_1}\})=0$, because of the fact that $x_{i_1}$ has no
partner. 
Therefore,   $d(y_{i_2})\leq a+1$ and hence, $y_{i_2}$ also has no partner. Continuing this process, as long as possible, as a result we obtain a path $P= y_{t+1}x_{i_1}y_{i_2}x_{i_2}\ldots x_{i_l}y_{i_l}$ or a cycle 
$C= y_{t+1}x_{i_1}y_{i_2}x_{i_2}\ldots x_{i_l}y_{t+1}$.  It is not difficult to see that all the vertices of this path (cycle) have no partner. If the former case holds, then $x_1$ is in $P$, which is a contradiction since $x_1$ has a partner (namely $x_2$).

If the second case holds, then, since every vertex of  $C$ has no partner in $D$, it follows that there is no arc from a vertex of $V(D)\setminus V(C)$ to a vertex of $V(C)$, which contradicts that $D$ is strongly connected. This completes the proof of the existence of a perfect matching from $X$ to $Y$. The proof for a perfect matching in the opposite direction is analogous. This completes the proof of the theorem.
\fbox \\\\

\section {Remarks}

Now using Theorem 1.10, we prove the following results (Lemmas 4.1-4.3).\\

\textbf{Lemma 4.1.} {\it Let $D$ be a strongly connected balanced bipartite digraph of order $2a\geq 8$ 
 with partite sets $X$ and $Y$ satisfying  condition $B_k$, $2\leq k\leq a/2$. If  $D$ is not Hamiltonian, then every vertex $u\in V(D)$ has a partner in $D$.} 

\textbf{Proof}. Let $D$ be a digraph  satisfying the conditions of the lemma. For a proof by contradiction, 
suppose that there is a vertex $x$ in $D$ which has no partner. By Theorem 1.10, $D$ has a cycle factor, say
 $C_1$, $C_2$, ... , $C_l$. Then $l\geq 2$ since $D$ is not Hamiltonian.  
Without loss of generality, we assume that $x\in V(C_1)$. It follows that $d^-(x^+_{C_1})=1$. Therefore,
 $d(x^+_{C_1})\leq a+1$. By Proposition 1,  this means that the vertex $x^+_{C_1}$ also has no partner. Similarly, we obtain that $d(x^{++}_{C_1})\leq a+1$ (where $x^{++}_{C_1}$ denotes the successor of $x^{+}_{C_1}$ on $C_1$) and hence, $x^{++}_{C_1}$ also has no partner in $D$.
Continuing this process,  we conclude that every vertex of $C_1$ has no partner in $D$. This  implies that there is no arc from a vertex of $A(V(D)\setminus V(C_1)$ to a vertex of $V(C_1))$,
 which contradicts that $D$ is strongly connected. The lemma  is proved. \fbox \\\\

\textbf{Lemma 4.2.}. {\it Let $D$ be a strongly connected balanced bipartite digraph of order $2a\geq 8$ 
 with partite sets $X$ and $Y$ satisfying  condition $B_k$, $2\leq k\leq a/2$. If  $D$ is not a cycle, then $D$ contains a non-Hamiltonian cycle of length at least four.}

\textbf{Proof}. Let $D$ be a digraph  satisfying the conditions of the lemma. For a proof by contradiction, 
suppose that $D$ contains non-Hamiltonian cycle of length at least four.

If $D$ is Hamiltonian, then it is not difficult to show that $D$ contains a non-Hamiltonian cycle of length at least 4. 
So we suppose, from now on, that $D$ is not Hamiltonian and contains no cycle of length at least 4. By Theorem 1.10, $D$ contains a cycle factor. Let $C_1, C_2, \dots, C_t$ be a minimal cycle factor of $D$ (i.e., $t$ is as small  as possible). Then the length of every $C_i$ is equal to two and $t=a$. Let $C_i=x_iy_ix_i$, where $x_i\in X$ and $y_i\in Y$. 
By Lemma 4.1, every vertex of $D$ has a partner. This means that  for every vertex 
$x\in V(D)$, $d(x)\geq a+k$ and $d^-(x)\geq k\geq 2$, $d^+(x)\geq k\geq 2$. Without loss of generality, we assume that 
$\{x_1,x_j\}$ with $j\not= 1$ is a dominating pair and $d(x_1)\geq 2a-k$. 

Let $Z$ be the subset of $Y$ with the maximum cardinality, such that  every vertex of $Z$ together with $x_1$ forms a cycle of length two. Without loss of generality, we assume that $Z=\{y_1,y_2,\ldots , y_l\}$. Then  $2a-k\leq d(x_1)\leq 2l+a-l=a+l$. Hence, $l\geq a-k$. Since $D$ contains no cycle of length four, it follows that the vertices $y_1$ and $x_i$, $2\leq i\leq l$, are not adjacent. Therefore, 
$$
a+k\leq d(y_1)\leq 2a-2l+2\leq 2k+2,
$$
i.e., $k\geq a-2$. Since $a/2\geq k\geq a-2$, we have   $a\geq 2k\geq 2a-4$, $a\leq 4$. 
If $a=4$, then $k=a/2=2$ and $l=a-k=2$. It is easy to see that $d(x_1)=d(y_1)=6$, the vertices $y_1$ and $x_i$, $3\leq i\leq 4$, form a cycle of length two and $x_1y_3\in A(D)$ or $y_3x_1\in A(D)$.
 Now it easy to see that $D$ contains a cycle of length four. Lemma 4.2 is proved. \fbox \\\\

For the next lemma we need the following lemma due to Bondy.\\

 \textbf{Bypass Lemma} (Lemma 3.17, Bondy \cite{[20]}). {\it Let $D$ be a strong   non-separable} ({\it i.e., $UG(D)$ is 2-connected}) {\it digraph, and let $H$ be a non-trivial proper subdigraph of $D$. Then $D$ contains an $H$-bypass}. \\

{\it Remark.} One can prove Bypass Lemma using the proof of Theorem 5.4.2 \cite{[1]}.\\

Now we will prove the following lemma.\\

\textbf{Lemma 4.3.} {\it Let $D$ be a strongly connected balanced bipartite digraph of order $2a\geq 8$ with partite sets $X$ and $Y$ satisfying condition $B_k$, where $2\leq k\leq a/2$. Then the following statements hold}:

(i){\it  the underlying undirected graph $UG(D)$ is 2-connected};

(ii) {\it if $C$ is a cycle of length $m$, $2\leq m\leq 2a-2$, then $D$ contains a $C$-bypass}.

\textbf{Proof}. (i). Suppose, on the contrary, that $D$ is a strongly connected balanced bipartite digraph of order $2a\geq 8$ with partite sets $X$ and $Y$ satisfies condition $B_k$ but $UG(D)$ is not 2-connected. Then $V(D)= E\cup F\cup \{u\}$, where $E$ and $F$ are non-empty  subsets, $E\cap F=\emptyset$, $u\notin E\cup F$, and there is no arc between $E$ and $F$. Since $D$ is strong, it follows that there are  vertices $x\in E$ and $y\in F$ such that $\{x,y\}\rightarrow u$, i.e., $\{x,y\}$ is a dominating pair. Without loss of generality, we assume that $x,y\in X$. Then $u\in Y$.
By condition $B_k$, it is easy to see that
$$
3a\leq d(x)+d(y)\leq 4+2|E\cap Y|+ 2|F\cap Y|\leq 2a+2,
$$
 which is a contradiction.  This proves that $UG(D)$ is 2-connected.

(ii).The second claim of the lemma is an immediate consequence of the first claim and  Bypass Lemma. Lemma 4.3 is proved. \fbox \\\\

Submitted 20.09.2017, accepted 16.01.2018.
\end{document}